\documentclass[11pt]{article}
\usepackage[a4paper, margin=1in]{geometry}
\usepackage[utf8]{inputenc}
\usepackage{amsmath}
\usepackage{amssymb}
\usepackage{amsthm}
\usepackage{mathtools}
\usepackage{graphicx}
\usepackage{enumitem}
\usepackage{float}
\newtheorem{theorem}{Theorem}[section]
\newtheorem{lemma}[theorem]{Lemma}
\newtheorem{corollary}[theorem]{Corollary}
\newtheorem{proposition}[theorem]{Proposition}
\newtheorem{claim}[theorem]{Claim}

\begin{document}
\title{Positive codegree thresholds for perfect matchings in hypergraphs}
\author{Richard Mycroft\thanks{School of Mathematics, University of Birmingham, UK. {\tt r.mycroft@bham.ac.uk}. RM is grateful for financial support from EPSRC Standard Grant EP/R034389/1.} \and
 Camila Z\'arate-Guer\'en\thanks{School of Mathematics, University of Birmingham, UK. {\tt ciz230@student.bham.ac.uk}.}}
\date{}
\maketitle

\begin{abstract}
We give, for each $k \geq 3$, the precise best possible minimum positive codegree condition for a perfect matching in a large $k$-uniform hypergraph $H$ on $n$ vertices. Specifically we show that, if $n$ is sufficiently large and divisible by $k$, and $H$ has minimum positive codegree $\delta^+(H) \geq \frac{k-1}{k}n - (k-2)$ and no isolated vertices, then $H$ contains a perfect matching. For $k=3$ this was previously established by Halfpap and Magnan, who also gave bounds for $k \geq 4$ which were tight up to an additive constant.
\end{abstract}

\section{Introduction}

A perfect matching in a graph (respectively $k$-uniform hypergraph, or $k$-graph for short) is a set of pairwise disjoint edges which collectively cover every vertex. So for any collection of objects -- which we represent by vertices -- and any property which can be mutually held by pairs ($k$-tuples) of objects -- which we represent by edges -- a perfect matching is precisely a partition of the objects into pairs ($k$-tuples) each of which holds the mutual property. In this way, the problem of finding a perfect matching in a graph or $k$-graph gives a versatile abstract framework which can be applied to an incredibly diverse range of scenarios in applications both within and outside mathematics. Consequently, this problem has been one of the most studied topics of research in extremal graph theory.

For graphs, celebrated theorems by Hall~\cite{Hall} and Tutte~\cite{Tutte} give elegant characterisations of graphs which contain perfect matchings, whilst Edmonds~\cite{Edmonds} gave an efficient algorithm to return such a matching or report that none exists. By contrast, for $k \geq 3$ the problem of finding a perfect matching in a $k$-graph was one of Karp's celebrated NP-hard problems~\cite{Karp}, and no short characterisation of $k$-graphs which contain perfect matchings is expected to exist. Consequently, most attention has focused on conditions which are sufficient for the existence of a perfect matching in a $k$-graph $H$. Since a necessary condition is that $H$ has no isolated vertices, the most natural types of condition to study are minimum degree condition. In this way we aim to generalise for $k$-graphs the elementary result that, for even $n$, every graph on $n$ vertices with minimum degree $\delta(G) \geq n/2$ contains a perfect matching.

The degree of a set $S \subseteq V(H)$ in a $k$-graph $H$, denoted $\deg(S)$, is the number of edges of~$H$ which contain $S$ as a subset. The minimum codegree of $H$, denoted $\delta(H)$, is then defined to be the minimum of $\deg(S)$ over all sets of $k-1$ vertices of $H$. This is the strongest commonly studied type of degree condition for $k$-graphs. The question of what minimum codegree forces a $k$-graph $H$ to contain a perfect matching was answered by R\"odl, Ruci\'nski and Szemer\'edi~\cite{RRS}, who improved on a previous asymptotically correct bound by K\"uhn and Osthus~\cite{KO} by giving the following theorem; a simple construction shows that this minimum codegree condition is best possible in a precise sense. Note also that the condition that $k$ divides $n$ is necessary since a matching of size $m$ in a $k$-graph covers precisely $mk$ vertices; we assume this condition without further comment for the rest of this discussion.

\begin{theorem}[{\cite[Theorem 1.1]{RRS}}]
For all $k\geq 3$ there exists $n_0$ for which the following holds for every $n \geq n_0$ which is divisible by $k$. If $H$ is a $k$-graph on $n$ vertices with 
$$\delta(H) \geq 
\begin{cases} 
\frac{n}{2} + 3 - k & \text{if } \frac{k}{2} \text{ is even and } \frac{n}{k} \text{ is odd}, \\
\frac{n}{2} + \frac{5}{2} - k & \text{if } k \text{ is odd and } \frac{n - 1}{2} \text{ is odd}, \\
\frac{n}{2} + \frac{3}{2} - k & \text{if } k \text{ is odd and } \frac{n - 1}{2} \text{ is even}, \\
\frac{n}{2} + 2 - k & \text{otherwise}.
\end{cases}$$ 
then $H$ contains a perfect matching.
\end{theorem}

One severe drawback of minimum codegree conditions is that they are simply too strong for many applications. Indeed, if a $k$-graph $H$ contains even a single pair of vertices $u, v$ for which no edge contains both $u$ and $v$, then $\delta(H) = 0$. This motivates the study of weaker notions of minimum degree. In particular, minimum $d$-degree conditions, which give a lower bound on the minimum of $\deg(S)$ over all sets of size $d$ for some $1 \leq d \leq k-2$, have been widely studied (see the surveys by K\"uhn and Osthus~\cite{KO2}, R\"odl and Ruci\'nski~\cite{RR} and Zhao~\cite{Zhao} for further details). However, these conditions lack the constructive power of minimum codegree conditions. Instead, Halfpap and Magnan~\cite{HM} considered perfect matchings of $k$-graphs with high minimum \emph{positive} codegree. The minimum positive codegree of a $k$-graph $H$, denoted $\delta^+(H)$ and introduced by Balogh, Lemons and Palmer~\cite{BLP}, is the minimum of $\deg(S)$ over all sets $S$ of $k$ vertices of $H$ with $\deg(S) \geq 1$ (that is, which are contained in some edge of $H$). By not placing any requirement on the degrees of sets which are not contained in an edge, minimum positive codegree conditions give a more versatile notion of minimum degree (for example, multipartite $k$-graphs may have non-trivial minimum positive codegree, but must have minimum codegree zero) whilst retaining similar constructive properties as minimum codegree conditions. Note, however, that minimum positive codegree conditions do not preclude the existence of isolated vertices, so we will also require the (very weak) condition that no vertex is isolated.

Halfpap and Magnan~\cite{HM} demonstrated that every $3$-graph on $n$ vertices with $\delta^+(H) \geq 2n/3 - 1$ and no isolated vertices contains a perfect matching, and that this minimum positive codegree condition is precisely best possible. For $k \geq 4$ they proved that, for large $n$, every $k$-graph on $n$ vertices with $\delta^+(H) \geq \frac{k-1}{k}n + k^2$ and no isolated vertices contains a perfect matching, and that this minimum positive codegree condition is best possible up to the constant additive error. The contribution of this paper is to determine, for large $n$, the exact best possible minimum positive codegree condition which ensures a perfect matching in a $k$-graph $H$ for all~$k$.

\begin{theorem} \label{main}
For all $k\geq 3$ there exists $n_0$ for which the following holds for every $n \geq n_0$ which is divisible by $k$. If $H$ is a $k$-graph on $n$ vertices with $\delta^+(H) \geq \frac{k-1}{k}n - (k-2)$ and with no isolated vertices then $H$ contains a perfect matching.
\end{theorem}

 The minimum positive codegree condition of Theorem~\ref{main} is best possible in the sense that for every $n \in \mathbb{N}$ which is divisible by $k$ there exists a $k$-graph on $n$ vertices with $\delta^+(H) = \frac{k-1}{k}n - (k-1)$ and with no isolated vertices which does not contain a perfect matching. This is shown by the following construction given by Halfpap and Magnan~\cite{HM}: let $A$ and $B$ be vertex-disjoint sets of size $|A| = \frac{n}{k} + 1$ and $|B| = n-|A|$, and let $H_\mathrm{ext}$ be the $k$-graph with vertex set $A \cup B$ (so $H_\mathrm{ext}$ has $n$ vertices) whose edges are all sets $e \in \binom{A \cup B}{k}$ with $|e \cap A| \in \{0,1\}$. We then have $\delta^+(H_\mathrm{ext}) = |B| - (k-2) = \frac{k-1}{k}n - (k-1)$, but $H$ has no perfect matching. Indeed, every edge of $H$ contains at most one vertex of $A$, so a matching in $H$ covering every vertex of~$A$ would contain at least $|A| > \frac{n}{k}$ edges, which is not possible.

While we were finalising this manuscript, Letzter and Ranganathan announced a proof that every $k$-graph on $n$ vertices with $\delta^+(H) \geq \frac{k-1}{k}n - (k-2)$ and no isolated vertices contains a tight Hamilton cycle, confirming a conjecture of Illingworth, Lang, M\"uyesser, Parczyk and Sgueglia~\cite{ILMPS}. This result would imply Theorem~\ref{main} since if $k$ divides $n$ then a tight cycle on $n$ vertices contains a perfect matching. However, we understand that their proof is likely to be substantially longer and more technical than the brief argument we present here.

 Our proof of Theorem~\ref{main} considers separately those graphs which are close to the construction $H_\mathrm{ext}$ and those which are far from $H_\mathrm{ext}$. With this in mind, we say that a $k$-graph $H$ on~$n$ vertices is \emph{$\gamma$-extremal} if there exists a set $S \subseteq V(H)$ of size $|S| = \frac{n}{k}$ for which at most $\gamma n^k$ edges of $H$ have at least two vertices in $S$.

\subsection{Notation and preliminaries}

For a set $X$ and $k \in \mathbb{N}$ we write $\binom{X}{k}$ for the set of all $k$-element subsets of $X$. When we assume $x \ll y$ we mean that for all $y > 0$ there exists $x_0$ such that for every $x > 0$ with $x \leq x_0$ the subsequent statements hold for $x$ and $y$, and the meaning of similar assumptions with more variables is analogous.

For a $k$-graph $H$, a set $E' \subseteq E(H)$ of edges of $H$ and a set $S \subseteq V(H)$, we write $\deg_{E'}(S)$ for the number of edges $e \in E'$ with $S \subseteq e$. So in particular $\deg_{E(H)}(S)$ is identical to $\deg_{H}(S)$.
We frequently use without further comment the fact that if $S \in \binom{V(H)}{k}$ then $\deg(S) = 1$ if $S \in E(H)$ and $\deg(S) = 0$ otherwise, so the statement $\deg(S) \geq 1$ means precisely that $S$ is an edge of $H$. 

We mostly apply minimum positive codegree conditions through the following proposition, which enables us to form edges one vertex at a time. This is a sharpened version of a similar proposition we gave in a previous work~\cite{MZG}; we need the sharp version to work with exact minimum positive codegree conditions.

\begin{proposition} \label{choosevs}
The following properties hold for every $k$-graph $H$.
\begin{enumerate}[label=(\roman*)]
    \item\label{choose1} For every $S \subseteq V(H)$ with $|S| \leq k-1$ and $\deg(S) \geq 1$ there are at least $\delta^+(H) + k - 1 - |S|$ vertices $x \in V(H)$ such that $\deg(S \cup \{x\}) \geq 1$.
    \item\label{choose2} Every vertex $v \in V(H)$ with $\deg(v)\geq 1$ has $\deg(v) \geq \binom{\delta^+(H) + k-2}{k-1}$.
\end{enumerate}
\end{proposition}
\begin{proof}
    For~\ref{choose1}, consider an edge $e$ containing $S$, which exists since $\deg(S) \geq 1$. Let $y\in e\setminus S $ and consider $T=e\setminus\{y\}$. We then have $|T \setminus S| = k-1-|S|$, and for each $x \in T \setminus S$ we have $S \cup \{x\} \subseteq T$, so $\deg(S\cup\{x\})\geq 1$.
    At the same time, there are at least $\delta^+(H)$ vertices $x\in V(H)\setminus T$ for which $T\cup\{x\} \in E(H)$ and so $\deg(S\cup\{x\})\geq 1$. Taking the union of both sets gives us the result.  
    
    For~\ref{choose2}, observe that repeated application of~\ref{choose1} gives us at least $(\delta^+(H)+k-2)!/(\delta^+(H)-1)!$ ordered sequences $(u_1, \dots, u_{k-1})$ such that $\{v, u_1, \dots, u_{k-1}\}$ is an edge of $H$. Since there are $(k-1)!$ permutations of a set $\{u_1, \dots, u_{k-1}\}$, the claimed bound follows. 
\end{proof}

We also use the following concentration bound for binomial distributions. 

\begin{theorem}[{\cite[Corollary 2.3]{JLR00:randomgraphs}}] \label{chernoff}
  If~$0<a<3/2$ and~$X$ is a binomial random variable then~$\mathbb{P}\bigl(\,|X-\mathbb{E} X|\geq a\mathbb{E} X\,\bigr)\leq 2\exp(-a^2 \mathbb{E} X/3)$.
\end{theorem}

\section{The extremal case}

We begin by showing that Theorem~\ref{main} holds for all $k$-graphs which are $\gamma$-extremal for sufficiently small $\gamma$. The argument for this proceeds by deleting a small matching $M$ in $H$ which covers all atypical vertices, then applying the following theorem of Daykin and H\"aggkvist~\cite{DH} to find a perfect matching in the subgraph of $H$ which remains after deleting the edges of $M$.

\begin{theorem}[{\cite[Theorem~1]{DH}}] \label{kpartitematching}
Let $H$ be a $k$-partite $k$-graph whose vertex classes each have size $n$. If every vertex $v$ of $H$ has $\deg(v) \geq (k-1)n^{k-1}/k$, then $H$ admits a perfect matching.
\end{theorem}

We also use the fact that every graph $G$ on $n$ vertices contains a matching covering at least $\min(n, 2 \delta(G))$ vertices, which is a standard exercise in elementary graph theory and can be proved by considering a maximum matching in $G$. Our use of this result in the proof of Lemma~\ref{extremal}, to obtain a matching $M_G$ in an auxiliary graph $G$, is the only point in the proof of Theorem~\ref{main} where we use the precise value of the minimum positive codegree condition.

\begin{lemma} \label{extremal}
For each $k \geq 3$ there exists $\gamma \in (0,1)$ for which the following holds for every $n > 8k$ which is divisible by $k$. If $H$ is a $k$-graph on $n$ vertices with $\delta^+(H) \geq \frac{k-1}{k}n - (k-2)$ and with no isolated vertices which is $\gamma$-extremal, then $H$ contains a perfect matching.     
\end{lemma}

\begin{proof}
Fix $k \geq 3$ and set $\gamma = 1/(2k)^{(2k)}$. Let $n > 8k$ be divisible by $k$, and let $H$ be a $k$-graph on $n$ vertices with $\delta^+(H) \geq \frac{k-1}{k}n - (k-2)$ which has no isolated vertices. Assume that $H$ is $\gamma$-extremal, so there is a set $S \subseteq V(H)$ with $|S| = \frac{n}{k}$ for which at most $\gamma n^k$ edges of $H$ have at least two vertices in $S$. Let $T = V(H) \setminus S$, let $K(S, T^{k-1})$ be the set of all sets $f \in \binom{V(H)}{k}$ with $|f \cap S| = 1$ and $|f \cap T| = k-1$, and let $F := K(S, T^{k-1}) \setminus E(H)$. Also let $E_{\geq 2}$ be the set of all edges of~$H$ with at least two vertices in $S$, and $E_1$ be the set of all edges of $H$ with precisely one vertex in $S$. Our choice of $S$ implies that $|E_{\geq 2}| \leq \gamma n^k$, so by Proposition~\ref{choosevs} we have 
\begin{align*}
    \sum_{v \in S} \deg_{E_1}(v) &= \sum_{v \in S} \deg_{H}(v) - \sum_{v \in S} \deg_{E_{\geq 2}}(v) \\ & \geq |S| \cdot \binom{\delta^+(H) + k-2}{k-1} - k |E_{\geq 2}| \geq |S| \cdot \binom{|T|}{k-1} - k \gamma n^k.
\end{align*}
It follows that $|F| \leq k \gamma n^k$, so at most $k^2 \sqrt{\gamma} n$ vertices $v \in V(H)$ have $\deg_F(v) > \sqrt{\gamma} n^{k-1}$.
Let $X \subseteq S$ and $Y \subseteq T$ consist of all vertices of $S$ and $T$  (respectively) with degree more than $\sqrt{\gamma} n^{k-1}$ in $F$. So $|X| + |Y| \leq k^2 \sqrt{\gamma} n$. Set $A = S \cup Y$ and $B = T \setminus Y$, so $|A| = \frac{n}{k}+|Y|$ and $|B| = \frac{k-1}{k}n-|Y|$.

Let $G$ be the graph with vertex set $A$ whose edges are all pairs $uv \in \binom{A}{2}$ with $\deg_H(uv) \geq 1$. By Proposition~\ref{choosevs} every $x \in A$ has $\deg_G(x) \geq \delta^+(H) + (k-2) - (n-|A|) \geq |Y|$. Since $|Y| \leq k^2\sqrt{\gamma} n \leq \frac{n}{2k} \leq \frac{|A|}{2}$ it follows that $G$ contains a matching $M_G$ of size $|Y|$. Let $L_G = (X \cup Y) \setminus V(M_G)$. Observe that for each pair $uv \in M_G$, since $\deg_H(uv) \geq 1$ by definition of $G$, we may extend $uv$ to an edge $e_{uv} \in E(H)$ with $|e_{uv} \cap B| = k-2$ by selecting vertices $b_1, b_2, \dots, b_{k-2} \in B$ in turn so that $\deg(uvb_1b_2\dots b_i) \geq 1$ for each $i \in [k-2]$. Indeed, by Proposition~\ref{choosevs} the number of options for the choice of each vertex is at least $\delta^+(H) - (n-|B|) > n/4$. In the same way, for each vertex $x \in L_G$ we may extend $\{x\}$ to an edge $e_{x} \in E(H)$ with $|e_{x} \cap B| = k-1$ by selecting vertices $b_1, b_2, \dots, b_{k-1} \in B$ in turn so that $\deg(xb_1b_2\dots b_i) \geq 1$ for each $i \in [k-1]$. Again Proposition~\ref{choosevs} implies that there are more than $n/4$ options for each choice. Since the edges $e_{uv}$ for $uv \in M_G$ and $e_x$ for $x \in L_G$ collectively cover at most $k(|M_G|+|L_G|) \leq k(|X|+2|Y|) \leq 2k^3 \sqrt{\gamma} n \leq n/4$ vertices in total, we may make these choices so that the edges in $M := \{e_{uv} : uv \in M_G\} \cup \{e_x : x \in L_G\}$ are pairwise vertex-disjoint, meaning that $M$ is a matching in $H$ of size $|Y|+|L_G|$ which covers every vertex of $X \cup Y$.

 Let $A' := A \setminus V(M)$, $B' := B \setminus V(M)$ and $H' = H[A' \cup B']$. So $H'$ has $n' := n - |V(M)| \geq n - 2k^3 \sqrt{\gamma} n \geq n/2$ vertices. Observe that also $n' = n - k(|Y|+ |L_G|)$, so $|A'| = |A| - 2|M_G| - |L_G| = \frac{n}{k} - (|Y| + |L_G|) = \frac{n'}{k}$ and $|B'| = |B| - (k-2)|M_G| - (k-1)|L_G| = \frac{k-1}{k}n - (k-1)(|Y|+|L_G|) = \frac{k-1}{k}n'$.
Arbitrarily partition $B'$ into $k-1$ sets of equal size, and let $H^*$ be the subgraph of $H$ with vertex set $A' \cup B'$ whose edges are all edges of $H$ with precisely one vertex in each part of~$B'$ and one vertex in $A'$. So $H^*$ is a $k$-partite $k$-graph whose vertex classes each have size $n'/k$. Moreover, since $A' \subseteq S \setminus X$ and $B' \subseteq T \setminus Y$, every vertex $v \in A' \cup B'$ has $\deg_{H^*}(v)\geq (\frac{n'}{k})^{k-1} - \deg_F(v) \geq (\frac{n'}{k})^{k-1} - \frac{n^{k-1}}{(2k)^k} \geq \frac{k-1}{k} (\frac{n'}{k})^{k-1}$, so $H^*$ contains a perfect matching $M'$ by Theorem~\ref{kpartitematching}. This gives a perfect matching $M \cup M'$ in $H$. 
\end{proof}

\section{Absorbing}

We next show that a substantially weaker minimum positive codegree condition than that of Theorem~\ref{main} ensures the existence of a small `absorbing' matching in a $k$-graph $H$.

\begin{lemma} \label{absorbing}
Suppose that $1/n \ll \beta \ll \alpha, 1/k$, and let $H$ be a $k$-graph on $n$ vertices. If $\delta^+(H) \geq \frac{n}{2} + \alpha n$ and $H$ has no isolated vertices, then there exists a set $A \subseteq V(H)$ of size $|A| \leq \beta n$ with the following absorbing property: for every $S \subseteq V(H) \setminus A$ of size $|S| \leq \beta^2 n$ whose size is divisible by $k$ there is a perfect matching in $H[A \cup S]$. 
\end{lemma}

\begin{proof}
We say a set of $k(k-1)$ vertices $W$ \textit{absorbs} a set $T$ of $k$ vertices  if both $H[W]$ and $H[W \cup T]$ contain perfect matchings (see Figure~\ref{absfig}). 
 \begin{figure}
\begin{center}
\includegraphics[width=3cm]{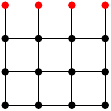}
\end{center} 
\caption{For $k=4$ the absorbing sets we form in Claim~\ref{absclaim} have the form shown. Observe that the set of 12 black vertices absorbs the set of 4 red vertices (the vertices in the top row).}
\label{absfig}
\end{figure}
We first prove that there are many absorbing sets for each set of $k$ vertices. For each $T \in \binom{V(H)}{k}$ let $\mathcal{A}(T)$ be the family of all sets $W \in \binom{V(H)}{k(k-1)}$ which absorb $T$.

\begin{claim} \label{absclaim}
    For each $T\in \binom{V(H)}{k}$ we have $|\mathcal{A}(T)|\geq (\alpha n)^{k(k-1)}/(k(k-1))!$.
\end{claim}

To prove the claim, let $T \in \binom{V(H)}{k}$ and write $T = \{t_1, \dots, t_k\}$. We construct pairwise disjoint edges $e_i = \{v_1^i, \dots, v_k^i\} \subseteq V(H)\setminus T$ for each $1\leq i\leq k-1$ in turn by choosing the vertices $v_j^i$ one by one. Specifically, for each $1\leq j\leq k$ we choose $v_j^i$ so that both $\deg(\{v^i_1, \dots, v^i_j\}) \geq 1$ and $ \deg(\{t_j, v_j^1, \dots, v^i_j\}) \geq 1$. Excluding the $k$ vertices of $T$ and the at most $k(k-1)$ vertices that we have previously chosen, by Proposition~\ref{choosevs} we then have at least $2\delta^+(H) - n - k^2 \geq \alpha n$ options for $v_j^i$, so we may indeed choose the vertices $v_j^i$ in this way. Let $W = \{v^i_j : i \in [k-1], j \in [k]\}$, so the edges $e_1, \dots, e_{k-1}$ form a perfect matching in $H[W]$.

Observe that our choices of the vertices $v_j^i$ also ensure that the sets $f_j := \{t_j, v_j^1, \dots, v_j^{k-1}\}$ for $j \in [k]$ are pairwise disjoint edges of $H$. So $\{f_1, \dots, f_k\}$ is a perfect matching in $H[W\cup T]$, and so $W \in \mathcal{A}(T)$. Since we had at least $\alpha n$ options for each of the $k(k-1)$ vertices in~$W$, we conclude that $|\mathcal{A}(T)| \geq (\alpha n)^{k(k-1)}/(k(k-1))!$. This completes the proof of Claim~\ref{absclaim}.
\medskip

Let $c := 3\beta^2 (k(k-1))!/\alpha^{k(k-1)}$. Choose a family $\mathcal{F} \subseteq \binom{V(H)}{k(k-1)}$ at random by including each $F \in \binom{V(H)}{k(k-1)}$ in $\mathcal{F}$ with probability $p = c n^{1-k(k-1)}$, independently of all other choices. The random variables $|\mathcal{F}|$ and $|\mathcal{A}(T)\cap \mathcal{F}|$ for each $T\in \binom{V(H)}{k}$ are then distributed binomially with expectations $\mathbb{E}[|\mathcal{F}|] \leq pn^{k(k-1)} = cn$ and $\mathbb{E}[|\mathcal{A}(T) \cap \mathcal{F}|] = p|\mathcal{A}(T)| \geq 3\beta^2 n$ respectively (the latter using Claim~\ref{absclaim}). So by Theorem~\ref{chernoff}, with high probability we have that $|\mathcal{F}| \leq 2cn \leq \beta n/k^2$ and that $|\mathcal{A}(T)\cap \mathcal{F}| \geq 2\beta^2 n$ for each set $T \in \binom{V(H)}{k}$.  

We denote by $I(\mathcal{F})$ the set of intersecting pairs in $\mathcal{F}$, that is, pairs $(F_1, F_2)\in \mathcal{F}\times \mathcal{F}$ such that $F_1\cap F_2 \neq \emptyset$. The random variable $|I(\mathcal{F})|$ has expectation $\mathbb{E}[{I(\mathcal{F})}] \leq p^2 \cdot (k(k-1))\cdot n^{k(k-1)}\cdot n^{k(k-1) - 1} = c^2k(k-1)n \leq \beta^2 n/2$. Hence by Markov's inequality we have $|I(\mathcal{F})| \leq \beta^2 n$ with probability at least $1/2$. So we may choose $\mathcal{F}$ for which all of these events hold, that is,
\begin{enumerate}[label = (\roman*), noitemsep]
    \item\label{F1} $|\mathcal{F}| \leq \beta n/k^2$,
    \item\label{F2} $|\mathcal{A}(T)\cap \mathcal{F}| \geq 2 \beta^2 n$ for each $T\in\binom{V(H)}{k}$, and
    \item\label{F3} there are at most $\beta^2 n$ intersecting pairs in $\mathcal{F}$. 
\end{enumerate}

Remove from $\mathcal{F}$ those sets which are not contained in $\mathcal{A}(T)$ for any $T\in\binom{V(H)}{k}$, and also one set from each pair in $I(\mathcal{F})$. Let $\mathcal{F}'$ be the resulting subset of $\mathcal{F}$, and let $A := \bigcup_{W\in \mathcal{F}'} V(W)$, so $|A| = |\mathcal{F}'| \cdot k(k-1) \leq k^2 |\mathcal{F}| \leq \beta n$ by~\ref{F1}. Now consider a set $S \subseteq V(H)\setminus A$ with $|S| = rk \leq \beta^2n$ for some $r\in\mathbb{N}$. Arbitrarily partition $S$ into $r$ sets $S_1, \dots, S_r$ of $k$ vertices each. By~\ref{F2} and~\ref{F3} we have $|\mathcal{A}(S_i)\cap \mathcal{F}'| \geq \beta^2 n$ for each $i\in[r]$, so
we may greedily choose a set $W_i \in \mathcal{A}(S_i)\cap \mathcal{F}'$ for each $i \in [r]$ so that the sets $W_1, \dots, W_r$ are all distinct. Finally, observe that for each $i\in[r]$ there is a perfect matching in $H[W_i\cup S_i]$ since $W_i \in \mathcal{A}(S_i)$, and for each $W \in \mathcal{F}' \setminus \{W_i : i \in [r]\}$ there is a perfect matching in $H[W]$ since our choice of $\mathcal{F}'$ ensures that $W$ absorbs some set. The union of these matchings is a perfect matching in $A \cup S$, completing the proof.
\end{proof}

\section{The non-extremal case}

Our focus in this section is on $k$-graphs $H$ which are not $\gamma$-extremal for some pre-fixed $\gamma$. Our aim is the following lemma, which states that if $H$ also satisfies a slightly weaker minimum positive codegree condition than that of Theorem~\ref{main}, then $H$ contains a matching covering almost all vertices of $H$.

\begin{lemma} \label{almostmatching}
Suppose that $1/n \ll \eta, \alpha \ll \gamma, 1/k$ and that $k$ divides $n$, and let $H$ be a $k$-graph on $n$ vertices. If $\delta^+(H) \geq \frac{k-1}{k}n - \alpha n$ and $H$ has no isolated vertices, then either $H$ is $\gamma$-extremal or $H$ contains a matching of size at least $(1-\eta) n/k$. 
\end{lemma}

Let $H$ be a $k$-graph. For vertices $u, v \in V(H)$ let $E_u(H)$ denote the set of edges of $H$ which contain $u$, and $E_{uv}(H)$ denote the set of edges of $H$ which contain both $u$ and $v$. A fractional matching in $H$ is a function $w: E(H) \to \mathbb{R}_{\geq 0}$ such that for every $v \in V(H)$ we have $\sum_{e \in E_v(H)} w(e) \leq 1$; $w$ is perfect if for every $v \in V(H)$ this inequality is actually an equality, which is equivalent to saying that $\sum_{e \in E(H)} w(e) = n/k$. For a pair of vertices $u,v\in E(H)$, we define $w_u := \sum_{e\in E_{u}} w(e)$ and $w_{uv} := \sum_{e\in E_{uv}} w(e)$. Our first step towards Lemma~\ref{almostmatching} is to show that the conditions of Lemma~\ref{almostmatching} are sufficient to ensure a perfect fractional matching in~$H$. To do this we use Farkas' lemma on solvability of finite systems of linear inequalities.

\begin{lemma}[Farkas' lemma]\label{farkas}
For every $\mathbf{v}\in\mathbb{R}^n$ and every finite set $\mathcal{X}\subseteq\mathbb{R}^n$, either there exist weights $\lambda_{\bf{x}} \geq 0$ for each $\bf{x} \in \mathcal{X}$ such that $\sum_{\bf{x} \in \mathcal{X}} \lambda_{\bf{x}} \bf{x} = \bf{v}$, or there exists $\mathbf{y}\in\mathbb{R}^n$ such that $\mathbf{y}\cdot\mathbf{x}\leq 0$ for every $\bf{x}\in\mathcal{X}$ and $\mathbf{y}\cdot\mathbf{v} > 0$.
\end{lemma}

\begin{lemma} \label{fracmatching} 
Suppose that $1/n \ll \alpha \ll \gamma, 1/k$ and that $k$ divides $n$, and let $H$ be a $k$-graph on $n$ vertices. If $\delta^+(H) \geq \frac{k-1}{k}n - \alpha n$ and $H$ has no isolated vertices, then either $H$ is $\gamma$-extremal or $H$ admits a perfect fractional matching. 
\end{lemma}

\begin{proof}
Let $\nu$ be such that $\alpha \leq \nu \leq \gamma/2$ and $\nu n$ is an integer. Enumerate the vertices of $V(H)$ as $v_1, \dots, v_n$, and for each set $S \in \binom{V(H)}{k}$ let $\mathbf{x}_S \in \mathbb{R}^n$ be the vector whose $i$th coordinate is $1$ if $v_i \in S$ and $0$ otherwise. Suppose that $H$ does not admit a perfect fractional matching. This means that there are no weights $w(e)\geq 0$ such that $\sum_{e\in E(H)} w(e)\cdot \mathbf{x}_e = \mathbf{1}$. Hence by Lemma~\ref{farkas} there exists $\mathbf{y} = \{y_1,\dots, y_n\}\in\mathbb{R}^n$ such that $\mathbf{y\cdot 1}>0$ and $\mathbf{y}\cdot \mathbf{x}_e \leq 0$ for all $e\in E(H)$.

By relabelling the vertices if necessary, we may assume without loss of generality that $y_i \leq y_j$ for all $i\leq j$. Write $L := \frac{k-1}{k}n$. For each $i \in [\nu n]$ let 
$$A_i : = \{v_{L+i}, v_{L-\nu n+i}, v_{(k-2)i-(k-3)}, v_{(k-2)i-(k-2)}, \dots, v_{(k-2)i}\}$$
and for each $i \in [\frac{n}{k} - \nu n]$ let
$$B_i := \{v_{n+1-i}, v_{(k-2)\nu n + (k-1)i-(k-2)}, v_{(k-2)\nu n + (k-1)i-(k-3)}, \dots, v_{(k-2)\nu n + (k-1)i}\}.$$
Observe that the sets $A_i$ and $B_i$ partition $V(H)$. Moreover, by repeated application of Proposition~\ref{choosevs} we may extend $\{v_n\}$ to an edge $e := \{v_n, v_{i_1}, \dots, v_{i_{k-1}}\}$ of $H$ with $i_j \geq \delta^+(H)$ for each~$j \in [k-1]$, from which it follows that $\mathbf{y} \cdot \mathbf{x}_e \geq \mathbf{y} \cdot \mathbf{x}_{B_i}$ for each $i \in [\frac{n}{k} - \nu n]$. If additionally there exists an edge $e' = \{v_{i'_1}, \dots, v_{i'_k}\}$ with $i'_1, i'_2 \geq L + \nu n$ and $i'_j \geq (k-2)\nu n$ for each $3 \leq j \leq k$, then we have $\mathbf{y} \cdot \mathbf{x}_{e'} \geq \mathbf{y} \cdot \mathbf{x}_{A_i}$ for each $i \in [\nu n]$, and so 
$$0 < \mathbf{y}\cdot\mathbf{1} = \mathbf{y}\cdot \Big(\sum_{i\in[\nu n]} \mathbf{x}_{A_i} + \sum_{i\in[n/k -\nu n]} \mathbf{x}_{B_i}\Big) \leq \nu n\cdot\mathbf{y}\cdot \mathbf{x}_{e'} + (n/k -\nu n)\cdot\mathbf{y}\cdot \mathbf{x}_e \leq 0,$$
where the first equality comes from the sets $A_i, B_i$ being a partition of $V(H)$. This gives a contradiction.

So we may assume that no such edge $e'$ exists. It follows that every pair $v_{i'_1}, v_{i'_2}$ with $i'_1, i'_2 \geq L + \nu n$ is in at most $(k-2) \nu n^{k-2}$ edges of $H$. Let $S = \{v_\ell : \ell > L\}$, so $|S| = n/k$. The number of edges of $H$ with at least two vertices in $S$ is then at most 
$$\binom{|S|}{2} (k-2) \nu n^{k-2} + \nu n \cdot n^{k-1} \leq 2\nu n^k,$$
so $H$ is $2\nu$-extremal and therefore $\gamma$-extremal, as required.
\end{proof}

We can in fact insist that for each pair $u, v$ of distinct vertices of $H$ the total weight of edges containing both $u$ and $v$ is small; this is the following corollary.  

\begin{corollary} \label{fracmatchingsmallweight}
Suppose that $1/n \ll \alpha, \varepsilon \ll \gamma, 1/k$ and that $k$ divides $n$, and let $H$ be a $k$-graph on~$n$ vertices. If $\delta^+(H) \geq \frac{k-1}{k}n - \alpha n$ and $H$ has no isolated vertices, then either $H$ is $\gamma$-extremal or~$H$ admits a perfect fractional matching in which $w_{uv} \leq 1/(\varepsilon n)$ for all distinct $u, v \in V(H)$.
\end{corollary}

\begin{proof}
 Let $W$ be the set of all perfect fractional matchings in $H$, and choose $w \in W$ which minimises $\max(w_{uv} : {uv \in \binom{V(H)}{2}})$; this minimum exists since $W$ is the feasible region of a linear program with integer coefficients. Let $M := \max_{u,v} w_{uv}$ be this minimum value, and fix $\nu<M$ so that for each $uv\in \binom{V(H)}{2}$ we have either $w_{uv} = M$ or $w_{uv} < M - \nu$.

Let $B \subseteq \binom{V(H)}{2}$ consist of all pairs $uv$ with $w_{uv} = M$, and let $H'$ be the spanning subgraph of $H$ formed by deleting from $H$ every edge $e\in E(H)$ for which some $uv \in B$ has $u, v \in e$. Suppose for a contradiction that $H'$ admits a perfect fractional matching~$w'$. We may then obtain a third perfect fractional matching $w''$ in $H$ by setting $w''(e) := (1-\nu)w(e) + \nu w'(e)$ for every $e\in E(H')$ and $w''(e) := (1-\nu)w(e)$ for each $e\in E(H) \setminus E(H')$. Each $uv \in \binom{V(H)}{2} \setminus B$ then has $w''_{uv} \leq (1-\nu)w_{uv} + \nu w'_{uv} < (1-\nu)(M-\nu) + \nu < M$, whilst each $uv \in B$ has $w''_{uv} = (1-\nu)w_{uv} < M$, giving a contradiction to our choice of $w$. We conclude that $H'$ does not admit a perfect fractional matching, so by Lemma~\ref{fracmatching} either $\delta^+(H') < \frac{k-1}{k}n - (\alpha + (k-1)\varepsilon) n$, or $H'$ has an isolated vertex, or $H'$ is $(\gamma-\varepsilon)$-extremal. 

If $\delta^+(H') < \frac{k-1}{k}n - (\alpha + (k-1)\varepsilon) n$, then some $S \in \binom{V(H)}{k-1}$ is contained in at least $\delta^+(H) - \delta^+(H') > (k-1)\varepsilon n$ edges in $E(H) \setminus E(H')$. Each such edge contains a pair in $B$, at least one vertex of which must be in $S$. We conclude that some vertex $u \in S$ is in at least $\varepsilon n$ pairs in $B$.

If instead $H'$ has an isolated vertex $v$, then $v$ is in at least $\deg_H(v) \geq (\frac{n}{2})^{k-1} \geq 2 \varepsilon n^{k-1}$ edges of $E(H) \setminus E(H')$, each of which contains a pair in $B$. For each $w$ with $vw \in B$ there are at most $n^{k-2}$ edges of $H$ containing $vw$, so either there are at least $\varepsilon n$ pairs in $B$ containing~$v$, or there are at least $\varepsilon n^{k-1}$ edges of $H$ which contain a pair $xy \in B$ with $x, y \neq v$. Since for each such $xy$ there are at most $n^{k-3}$ edges of $H$ containing $x, y$ and $v$, it follows in the latter case that $|B| \geq \varepsilon n^2$, so some vertex $u \in V(H)$ is in at least $\varepsilon n$ pairs in $B$. 

Finally, if $H'$ is $(\gamma-\varepsilon)$-extremal then this fact is witnessed by a set $S \subseteq V(H)$.  
If $|E(H) \setminus E(H')| < \varepsilon n^k$ then $S$ also witnesses that $H$ is $\gamma$-extremal. On the other hand,  if $|E(H) \setminus E(H')| \geq \varepsilon n^k$ then $|B| \geq \varepsilon n^2$, so some vertex $u \in V(H)$ is in at least $\varepsilon n$ pairs in $B$. 

We conclude that either $H$ is $\gamma$-extremal or there is a vertex $u \in V(H)$ contained in at least~$\varepsilon n$ pairs in $B$. The latter implies that $1 = \sum_{e \in E_u(H)} w(e) \geq M \varepsilon n$, so $M \leq 1/(\varepsilon n)$.
\end{proof} 

Our final tool is a well-known theorem by Pippenger and Spencer~\cite{PS} which states that every approximately regular multi-$k$-graph with small pair degrees contains a matching covering almost all vertices. To fit with our arguments on fractional matchings we use the following restatement for fractional matchings given by Bowtell, Kathapurkar, Morrison and Mycroft~\cite{BKMM}.

\begin{theorem}[{\cite[Corollary 4.5]{BKMM}}] \label{frac_to_almost} 
Fix $k \geq 2$, suppose that $1/n \ll \varepsilon \ll \eta, 1/k$, and let $H$ be a $k$-graph on $n$ vertices. If $H$ admits a fractional matching in which $w_u \geq 1-\varepsilon$ for every $u \in V(H)$ and $w_{uv} \leq \varepsilon$ for all distinct $u, v \in V(H)$, then $H$ contains a matching of size at least $(1-\eta) n/k$.
\end{theorem}

We now give the proof of the main lemma of this section.

\begin{proof}[Proof of Lemma~\ref{almostmatching}]
Introduce a new constant $\varepsilon$ with $1/n \ll \varepsilon \ll \eta, \alpha$. By Corollary~\ref{fracmatchingsmallweight} either $H$ is $\gamma$-extremal or $H$ admits a perfect fractional matching in which $w_{uv} \leq \varepsilon$ for all distinct $u, v \in V(H)$; in the latter case $H$ contains a matching of size at least $(1-\eta) n/k$ by Theorem~\ref{frac_to_almost}.
\end{proof}

\section{Proof of Theorem~\ref{main}}
Fix $k \geq 3$ and choose $\gamma \in (0,1)$ according to Lemma~\ref{extremal}. Next choose $\beta$, $n_{\ref{almostmatching}}$ and $n_{\ref{absorbing}}$ with $\beta \leq \gamma/4$ so that 
\begin{enumerate}[noitemsep, label=(\arabic*)]
    \item for every $n \geq n_{\ref{absorbing}}$ we may apply Lemma~\ref{absorbing} with this value of $\beta$ and with $\alpha = 1/10$, and
    \item for every $n \geq n_{\ref{almostmatching}}$ we may apply Lemma~\ref{almostmatching} with $\gamma/2$, $2\beta$ and $\beta^2$ in place of $\gamma$, $\alpha$ and $\eta$ respectively.
\end{enumerate}
Our applications of Lemmas~\ref{absorbing} and~\ref{almostmatching} in the rest of the argument will be with the constants set out in (1) and (2). Finally set $n_0 := 2\max\{n_{\ref{almostmatching}}, n_{\ref{absorbing}}, 5k, k/\beta\}$.

Fix $n \geq n_0$ which is divisible by $k$, and let $H$ be a $k$-graph on $n$ vertices with $\delta^+(H) \geq \frac{k-1}{k} n -(k-2)$ and with no isolated vertices. Apply Lemma~\ref{absorbing} to obtain a set $A \subseteq V(H)$ of size $|A| \leq \beta n$ such that for every $S \subseteq V(H) \setminus A$ of size $|S| \leq \beta^2 n$ whose size is divisible by $k$ there is a perfect matching in $H[A \cup S]$. In particular, taking $S = \emptyset$ we find that $H[A]$ contains a perfect matching, so $|A|$ is divisible by $k$. Let $H' := H[V(H) \setminus A]$ and $n' := n - |A| \geq (1-\beta) n \geq n_{\ref{almostmatching}}$, so $|V(H')| = n'$ is divisible by $k$ since the same is true of both $n$ and $|A|$. Observe that $H'$ has no isolated vertices and $\delta^+(H') \geq \delta^+(H) - |A| \geq \frac{k-1}{k} n - 2 \beta n \geq \frac{k-1}{k} n' - 2 \beta n'$. So by Lemma~\ref{almostmatching} either $H'$ is $\gamma/2$-extremal or $H'$ contains a matching $M'$ of size at least $(1-\beta^2)n'/k$. 

Suppose first that $H'$ contains a matching $M'$ of size at least $(1-\beta^2)n'/k$, and let $S$ be the set of vertices of $H'$ not covered by $M'$, so $|S| \leq \beta^2 n$. Moreover, since $|V(M')|$ and $|V(H')|$ are divisible by $k$, the same is true of $|S|$. So by our choice of $A$ there is a perfect matching $M$ in $H[A \cup S]$, and so we have a perfect matching $M \cup M'$ in $H$.

Now suppose instead that $H'$ is $\gamma/2$-extremal, meaning that there is a set $S' \subseteq V(H')$ of size $n'/k$ for which at most $(\gamma/2) (n')^k \leq \gamma n^k/2$ edges of $H'$ have at least two vertices in~$S'$. Arbitrarily choose a set $X$ of $(n-n')/k$ vertices of $V(H) \setminus S'$ and set $S = S' \cup X$, so $|S| = n/k$. Every edge $e \in E(H)$ with $|e \cap S| \geq 2$ must then either contain a vertex of $X \cup A$ or be one of the at most $\gamma n^k/2$ edges of $H'$ with at least two vertices in~$S'$. Since at most $(|X|+|A|) n^{k-1} \leq 2\beta n^k \leq \gamma n^k/2$ edges of $H$ contain a vertex of $X \cup A$, we conclude that at most $\gamma n^k$ edges of $H$ have at least two vertices in $S$. So the set  $S$ witnesses that $H$ is $\gamma$-extremal, and it follows by Lemma~\ref{extremal} that $H$ contains a perfect matching. \qed

\bibliographystyle{plain}
\bibliography{biblio2}
\end{document}